# An implicit integration factor method for a kind of spatial fractional diffusion equations

**Yong-Liang Zhao[1, 3], Pei-Yong Zhu[1], Xian-Ming Gu[2], Xi-Le Zhao[1] and Huan-Yan Jian[1]**

1. School of Mathematical Sciences, University of Electronic Science and Technology of China, Chengdu, Sichuan 611731, P.R. China
2. School of Economic Mathematics/Institute of Mathematics, Southwestern University of Finance and Economics, Chengdu, Sichuan 611130, P.R. China
3. Email: ylzhaofde@sina.com

**Abstract**. A kind of spatial fractional diffusion equations in this paper are studied. Firstly, an $L1$ formula is employed for the spatial discretization of the equations. Then, a second order scheme is derived based on the resulting semi-discrete ordinary differential system by using the implicit integration factor method, which is a class of efficient semi-implicit temporal scheme. Numerical results show that the proposed scheme is accurate even for the discontinuous coefficients.

## 1. Introduction

In this work, a special nonlinear case of the spatial fractional diffusion equation (SFDE) proposed in [1] is studied:

$$\begin{cases} \dfrac{\partial u(x,t)}{\partial t} + {}_a^C D_x^\alpha \left( d_+(x) {}_x D_b^\alpha u(x,t) \right) + {}_x^C D_b^\alpha \left( d_-(x) {}_a D_x^\alpha u(x,t) \right) = f\left(u(x,t)\right), \\ \qquad\qquad\qquad\qquad\qquad\qquad\qquad\qquad (x,t) \in [a,b] \times [0,T], \\ u(a,t) = u(b,t) = 0, \qquad\qquad\qquad\qquad\qquad 0 \leq t \leq T, \\ u(x,0) = u_0(x), \qquad\qquad\qquad\qquad\qquad\qquad a \leq x \leq b, \end{cases} \qquad (1)$$

where $\alpha \in (\dfrac{1}{2},1)$, $d_\pm(x) \geq 0$ and $f\left(u(x,t)\right)$ is the nonlinear reaction term. The space fractional derivatives ${}_a^C D_x^\alpha u(x,t)$ and ${}_x^C D_b^\alpha u(x,t)$ ($0 < \alpha < 1$) are introduced in the Caputo sense [2]:

$${}_a^C D_x^\alpha u(x,t) = \dfrac{1}{\Gamma(1-\alpha)} \int_a^x (x-\eta)^{-\alpha} \dfrac{\partial u(\eta,t)}{\partial \eta} d\eta, \quad {}_x^C D_b^\alpha u(x,t) = \dfrac{-1}{\Gamma(1-\alpha)} \int_x^b (\eta-x)^{-\alpha} \dfrac{\partial u(\eta,t)}{\partial \eta} d\eta,$$

where $\Gamma(\cdot)$ denotes the Gamma function. While ${}_a D_x^\alpha u(x,t)$ and ${}_x D_b^\alpha u(x,t)$ ($0 < \alpha < 1$) represent the left- and right-Riemann-Liouville fractional derivatives [2] defined as:

$${}_a D_x^\alpha u(x,t) = \dfrac{1}{\Gamma(1-\alpha)} \dfrac{d}{dx} \int_a^x \dfrac{u(\eta,t)}{(x-\eta)^\alpha} d\eta, \quad {}_x D_b^\alpha u(x,t) = \dfrac{-1}{\Gamma(1-\alpha)} \dfrac{d}{dx} \int_x^b \dfrac{u(\eta,t)}{(\eta-x)^\alpha} d\eta,$$

The differential equations with nonlocal derivatives known as fractional differential equations (FDEs) have been attracted many researchers' interest in recent decades. Due to the nonlocal feature of fractional derivatives, FDEs are more suitable than the traditional differential equations in the description of the anomalous diffusion phenomena. Nowadays, the applications of FDEs have been recognized in numerous fields such as the physics [3], American options pricing [4], entropy [5] and image processing [6]. It is noticeable that finding the closed-form analytical solutions of most FDEs poses a challenge for researchers. For this reason, abundant numerical methods, e.g., finite difference method [7-10], finite element method [11, 12], and meshless method [13], have been proposed to solve the FDEs. Moreover, numerous fast solvers are designed based on the structures of the resulting numerical schemes, readers are suggested to refer to [7, 14-20] and the references therein.

All the above mentioned studies are based on the fully discrete schemes. In this work, a semi-implicit scheme is developed to approximate (1). To our best knowledge, no article considers such the semi-implicit scheme to solve SFDE (1) by employing the implicit integration factor (IIF) method [21]. Traditional integrating factor [22] or exponential time differencing methods [23-25] are not efficient for the systems with severely stiff reactions, because they still treat the reaction terms explicitly. To address this problem, Nie et al. [21] proposed a new class of semi-implicit schemes for stiff systems. They showed that the stability region of their methods is much larger than existing methods, who treat the reactions explicitly. Furthermore, numerical results in [21] demonstrate that their proposed schemes are accurate, robust and efficient. Later, Nie et al. [26] proposed a compact IIF method to solve the high-dimensional stiff reaction-diffusion equations. Such the method preserves the stability property of the IIF method and saves the storage requirement and CPU times. Other studies about IIF method can be found in [27-31].

The rest of this paper is organized as follows. In Section 2, a second-order implicit integration factor (IIF2) scheme is derived and its linear stability is also studied. Several numerical examples are provided in Section 3. Some conclusions are drawn in Section 4.

## 2. A second-order implicit integration factor scheme
In this section, the IIF2 scheme is derived and its stability is studied.

*2.1. The second-order scheme*

Let $h=(b-a)/N$ be the spatial step size for the positive integer $N$. Then the space domain $[a,b]$ can be covered by $\omega_h = \{x_i \mid x_i = a+ih, 0 \leq i \leq N\}$. Using the $L1$ formula [32]:

$$_a^C D_x^\alpha u(x)|_{x=x_i} \approx \frac{h^{-\alpha}}{\Gamma(2-\alpha)}\left[ a_0^{(\alpha)} u(x_i) - \sum_{k=1}^{i-1}\left(a_{i-k-1}^{(\alpha)} - a_{i-k}^{(\alpha)}\right)u(x_k) - a_{i-1}^{(\alpha)} u(x_0) \right]$$

with $a_i^{(\alpha)} = (i+1)^{1-\alpha} - i^{1-\alpha}$ $(i=0,1,2,\cdots)$ and the relationships between Caputo and Riemann-Liouville fractional derivatives [2]:

$$_a D_x^\alpha u(x) = {_a^C D_x^\alpha u(x)} + \frac{(x-a)^{-\alpha}}{\Gamma(1-\alpha)} u(a), \quad _x D_b^\alpha u(x) = {_x^C D_b^\alpha u(x)} + \frac{(b-x)^{-\alpha}}{\Gamma(1-\alpha)} u(b),$$

the spatial derivatives at $x=x_i$ $(i=1,2,\cdots,N-1)$ are approximated respectively as [1]:

$$_a^C D_x^\alpha \left( d_+(x) {_x D_b^\alpha u(x,t)} \right)|_{x=x_i} \approx \eta \left( \sum_{s=0}^{i-1} g_s^{(\alpha)} d_{+,i-s} \sum_{k=0}^{N-1-i+s} g_k^{(\alpha)} u_{i-s+k}(t) - a_{i-1}^{(\alpha)} d_{+,0} \sum_{k=1}^{N-1} g_k^{(\alpha)} u_k(t) \right),$$

$$_x^C D_b^\alpha \left( d_-(x) {_a D_x^\alpha u(x,t)} \right)|_{x=x_i} \approx \eta \left( \sum_{s=0}^{N-1-i} g_s^{(\alpha)} d_{-,i+s} \sum_{k=0}^{i+s-1} g_k^{(\alpha)} u_{i+s-k}(t) - a_{N-1-i}^{(\alpha)} d_{-,N} \sum_{k=1}^{N-1} g_k^{(\alpha)} u_{N-k}(t) \right).$$

Here $u_i(t) \approx u(x_i,t)$, $\eta = 1/h^{2\alpha}$, $d_{\pm,i} = d_\pm(x_i)$, $a_i^{(\alpha)} = \frac{1}{\Gamma(2-\alpha)}\left[(i+1)^{1-\alpha} - i^{1-\alpha}\right]$ and

$$g_k^{(\alpha)} = \begin{cases} a_0^{(\alpha)}, & k = 0, \\ a_k^{(\alpha)} - a_{k-1}^{(\alpha)}, & k > 0. \end{cases}$$

Then the semi-discrete ordinary differential system of Eq. (1) is

$$\frac{du_i(t)}{dt} + \eta \left( \sum_{s=0}^{i-1} g_s^{(\alpha)} d_{+,i-s} \sum_{k=0}^{N-1-i+s} g_k^{(\alpha)} u_{i-s+k}(t) - a_{i-1}^{(\alpha)} d_{+,0} \sum_{k=1}^{N-1} g_k^{(\alpha)} u_k(t) \right.$$
$$\left. + \sum_{s=0}^{N-1-i} g_s^{(\alpha)} d_{-,i+s} \sum_{k=0}^{i+s-1} g_k^{(\alpha)} u_{i+s-k}(t) - a_{N-1-i}^{(\alpha)} d_{-,N} \sum_{k=1}^{N-1} g_k^{(\alpha)} u_{N-k}(t) \right) = f(u_i(t)),\ 1 \leq i \leq N-1. \quad (2)$$

Denoting

$$\mathbf{u}(t) = [u_1(t), u_2(t), \cdots, u_{N-1}(t)]^T,\quad \frac{d\mathbf{u}(t)}{dt} = \left[ \frac{du_1(t)}{dt}, \frac{du_2(t)}{dt}, \cdots, \frac{du_{N-1}(t)}{dt} \right]^T,$$

$$D_+ = diag(d_{+,0}, d_{+,1}, \cdots, d_{+,N-1}),\quad D_- = diag(d_{-,1}, d_{-,2}, \cdots, d_{-,N}),$$

$$\mathbf{a}_+^{(\alpha)} = -[a_0^{(\alpha)}, a_1^{(\alpha)}, \cdots, a_{N-2}^{(\alpha)}]^T,\quad \mathbf{a}_-^{(\alpha)} = -[a_{N-2}^{(\alpha)}, a_{N-3}^{(\alpha)}, \cdots, a_0^{(\alpha)}]^T,$$

$$\mathbf{g}_+^{(\alpha)} = [g_1^{(\alpha)}, g_2^{(\alpha)}, \cdots, g_{N-1}^{(\alpha)}]^T,\quad \mathbf{g}_-^{(\alpha)} = [g_{N-1}^{(\alpha)}, g_{N-2}^{(\alpha)}, \cdots, g_1^{(\alpha)}]^T.$$

Then let

$$G_{L_+} = [\mathbf{a}_+^{(\alpha)}\ \tilde{G}],\ G_{R_+} = [\mathbf{g}_+^{(\alpha)}\ \tilde{G}]^T,\ G_{L_-} = [\tilde{G}^T\ \mathbf{a}_-^{(\alpha)}]^T,\ G_{R_-} = \begin{bmatrix} \tilde{G} \\ (\mathbf{g}_-^{(\alpha)})^T \end{bmatrix},$$

in which

$$\tilde{G} = \begin{bmatrix} g_0^{(\alpha)} & 0 & \cdots & 0 \\ g_1^{(\alpha)} & g_0^{(\alpha)} & \ddots & \vdots \\ \vdots & \ddots & \ddots & 0 \\ g_{N-2}^{(\alpha)} & \cdots & g_1^{(\alpha)} & g_0^{(\alpha)} \end{bmatrix}.$$

With the help of these symbols, the matrix form of Eq. (2) can be written as:

$$\frac{d\mathbf{u}(t)}{dt} + A\mathbf{u}(t) = \mathbf{f}(\mathbf{u}(t)), \quad (3)$$

where $A = \eta(G_{L_+} D_+ G_{R_+} + G_{L_-} D_- G_{R_-})$ and $\mathbf{f}(\mathbf{u}(t)) = [f(u_1(t)), f(u_2(t)), \cdots, f(u_{N-1}(t))]^T$.

For the discretization of the time direction, we use the IIF method proposed in [21] instead of the finite difference method. Let $\tau = \frac{T}{M}$ for the positive integer $M$, then $t_j = j\tau\ (0 \leq j \leq M)$. According to the work [21], the IIF2 scheme of (3) is

$$\mathbf{u}^{j+1} = e^{-A\tau}\left( \mathbf{u}^j + \frac{\tau}{2}\mathbf{f}(\mathbf{u}^j) \right) + \frac{\tau}{2}\mathbf{f}(\mathbf{u}^{j+1}), \quad (4)$$

where $\mathbf{u}^j = \mathbf{u}(t_j)$ and $\mathbf{f}(\mathbf{u}^j) = \mathbf{f}(\mathbf{u}(t_j))$.

### 2.2. Linear stability analysis of IIF2

Similar to Section 3 in [21], the linear stability of the IIF2 scheme is studied in this subsection. Applying Eq. (4) to the following scalar linear equation:

$$u_t = -qu + ru \text{ with } q > 0,$$

where $q$ and $r$ are diffusion and reaction, respectively. Then substituting $\mathbf{u}^j = e^{\iota j\theta}$ ($\iota$ is the imaginary unit) into the resulting equation. After some simple manipulations, we have
$$e^{\iota\theta} = e^{-q\tau}(1+\lambda/2) + \lambda e^{\iota\theta}/2, \quad 0 \leq \theta \leq 2\pi, \tag{5}$$
where $\lambda = r\tau$ with its real part $\lambda_r$ and imaginary part $\lambda_t$. Substituting $\lambda = \lambda_r + \iota\lambda_t$ into Eq. (5) and noticing $e^{\iota\theta} = \cos\theta + \iota\sin\theta$, we obtain
$$\lambda_r = \frac{2(1-e^{-2q\tau})}{c(q\tau,\theta)}, \quad \lambda_t = \frac{4\sin\theta e^{-q\tau}}{c(q\tau,\theta)},$$
in which $c(q\tau,\theta) = (1-e^{-q\tau})^2 + 2(1+\cos\theta)e^{-q\tau}$.

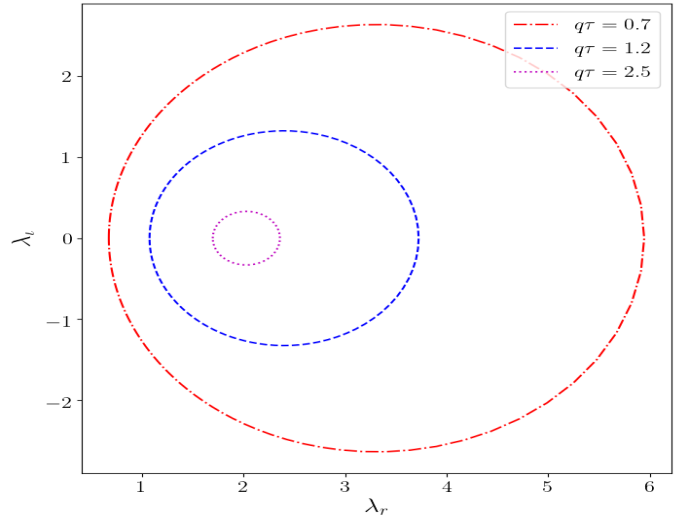

**Figure 1.** Stability regions (exterior of the closed curves) for (4).

Figure 1 is plotted to show the stability regions (exterior of the closed curves) for Eq. (4) with $q\tau = 0.7, 1.2, 2.5$. From such the figure, if $q \to 0$, the stability region becomes the domain $\lambda_r < 0$. If $q \to +\infty$, the stability region approaches the whole complex plane excluding the point $(2,0)$. Moreover, the IIF2 scheme (4) is $A$-stable because its stability region includes the domain $\lambda_r < 0$ in the complex plane of all $\lambda$.

## 3. Numerical experiments

In this section, several numerical examples are provided to show the order of accuracy of our scheme. Let $\xi_i^j = u(x_i,t_j) - u_i(t_j)$ be the error between the exact and numerical solutions of Eq. (1). Then denote
$$error(h,\tau) = \max_{0\leq i\leq N, 0\leq j\leq M}|\xi_i^j|, \quad rate_1 = \log_{\tau_1/\tau_2}\frac{error(h,\tau_1)}{error(h,\tau_2)}, \quad rate_2 = \log_{h_1/h_2}\frac{error(h_1,\tau)}{error(h_2,\tau)}.$$

In this work, the fixed-point iteration method (the maximum number of iterations and tolerance are respectively 200, 1e-12) is used to solve Eq. (4). All experiments were performed on a Windows 7 (32 bit) desktop-Intel(R) Core(TM) i3-2130 CPU 3.40GHz, 4GB of RAM using Spyder 3.2.8.

**Example 1** In this example, we consider the SFDE (1) with $a=0, b=1, T=1$, $d_+(x) = (1-x)^\alpha$, $d_-(x) = x^\alpha$, $u_0(x) = 10x^2(1-x)^2$ and $f(u) = 100u(u-0.5)(1-u)$. Since the exact solution of this case is unknown, the numerical solution computed from the mesh $M = N = 1024$ is treated as our exact solution.

**Table 1.** Maximum norm errors and convergence orders for Example 1 where $h = 1/1024$.

| $\tau$ | $\alpha = 0.6$ | | $\alpha = 0.7$ | | $\alpha = 0.9$ | |
|---|---|---|---|---|---|---|
| | $error(h,\tau)$ | $rate_1$ | $error(h,\tau)$ | $rate_1$ | $error(h,\tau)$ | $rate_1$ |
| 1/32 | 1.6739E-02 | - | 1.2248E-02 | - | 8.4310E-03 | - |
| 1/64 | 4.1883E-03 | 1.9988 | 3.4843E-03 | 1.8136 | 2.4637E-03 | 1.7749 |
| 1/128 | 1.1086E-03 | 1.9177 | 9.0807E-04 | 1.9400 | 6.3026E-04 | 1.9668 |
| 1/256 | 2.6787E-04 | 2.0491 | 2.2132E-04 | 2.0367 | 1.5643E-04 | 2.0104 |

**Table 2.** Maximum norm errors and convergence orders for Example 1 where $\tau = 1/1024$.

| $h$ | $\alpha = 0.6$ | | $\alpha = 0.7$ | | $\alpha = 0.9$ | |
|---|---|---|---|---|---|---|
| | $error(h,\tau)$ | $rate_2$ | $error(h,\tau)$ | $rate_2$ | $error(h,\tau)$ | $rate_2$ |
| 1/16 | 3.5064E-02 | - | 2.8995E-02 | - | 1.8068E-02 | - |
| 1/32 | 1.5180E-02 | 1.2079 | 1.3325E-02 | 1.1217 | 9.3471E-03 | 0.9509 |
| 1/64 | 5.8721E-03 | 1.3702 | 5.6120E-03 | 1.2476 | 4.5109E-03 | 1.0511 |
| 1/128 | 2.0759E-03 | 1.5002 | 2.2091E-03 | 1.3451 | 2.0478E-03 | 1.1394 |

It can be seen clearly from Table 1 that when $h = 1/1024$, the $error(h,\tau)$ decreases steadily with the shortening of $\tau$, and the order of accuracy in time is two. Fixing $\tau = 1/1024$, Table 2 lists the maximum norm errors and illustrates that the spatial convergence order is of $h^{2-\alpha}$. In a word, Tables 1-2 confirm that the rate of the truncation error of the IIF2 scheme (4) is $O(\tau^2 + h^{2-\alpha})$.

**Table 3.** Maximum norm errors and convergence orders for Example 2 where $h = 1/1024$.

| $\tau$ | $\alpha = 0.6$ | | $\alpha = 0.7$ | | $\alpha = 0.9$ | |
|---|---|---|---|---|---|---|
| | $error(h,\tau)$ | $rate_1$ | $error(h,\tau)$ | $rate_1$ | $error(h,\tau)$ | $rate_1$ |
| 1/64 | 3.3441E-05 | - | 4.9016E-05 | - | 9.4517E-05 | - |
| 1/128 | 9.7701E-06 | 1.7752 | 1.3286E-05 | 1.8834 | 3.1002E-05 | 1.6082 |
| 1/256 | 2.9065E-06 | 1.7491 | 3.5182E-06 | 1.9169 | 7.9265E-06 | 1.9676 |
| 1/512 | 7.3333E-07 | 1.9868 | 7.9535E-07 | 2.1452 | 1.6742E-06 | 2.2432 |

**Example 2** Considering the SFDE (1.1) with $a = -1, b = 1, T = 1$, discontinuous coefficients

$$d_+(x) = \begin{cases} 1.5e^{-x}, & -1 \leq x < 0, \\ 1, & 0 \leq x \leq 1, \end{cases} \quad d_-(x) = 1 \, (-1 \leq x \leq 1),$$

$u_0(x) = \dfrac{4e^{10x}}{\left(e^{10x}+1\right)^2}$ and $f(u) = \sin u$. For this case, no exact solution exists. Similar to Example 1, we also treat the numerical solution, which calculated from the mesh $M = N = 1024$, as the exact solution.

It can be seen from Tables 3-4 that the order of accuracy of the IIF2 scheme (4) is also $O(\tau^2 + h^{2-\alpha})$ for the discontinuous coefficients. As a conclusion, our scheme is accurate for the discontinuous coefficients.

Table 4. Maximum norm errors and convergence orders for Example 2 where $\tau = 1/1024$.

| h | $\alpha = 0.6$ | | $\alpha = 0.7$ | | $\alpha = 0.9$ | |
|---|---|---|---|---|---|---|
| | $error(h,\tau)$ | $rate_2$ | $error(h,\tau)$ | $rate_2$ | $error(h,\tau)$ | $rate_2$ |
| 1/64 | 9.0058E-03 | - | 8.3461E-03 | - | 4.1270E-03 | - |
| 1/128 | 3.9670E-03 | 1.1828 | 3.8174E-03 | 1.1285 | 2.1459E-03 | 0.9435 |
| 1/256 | 1.6806E-03 | 1.2391 | 1.5441E-03 | 1.3058 | 9.5260E-04 | 1.1716 |
| 1/512 | 5.5843E-04 | 1.5895 | 4.8088E-04 | 1.6830 | 3.1847E-04 | 1.5807 |

## 4. Concluding remarks

In this work, the IIF2 scheme (4) is derived by employing the $L1$ formula and the IIF method, which treats the diffusion term exactly and the nonlinear reaction term implicitly, to approximate a kind of nonlinear SFDE (1). Then, the linear stability of (4) is analysed in Section 2.2. Numerical results show that the order of accuracy of our scheme is of $O(\tau^2 + h^{2-\alpha})$ for the continuous/discontinuous coefficients.

### Acknowledgments

This research is supported by the National Natural Science Foundation of China (Nos. 61772003, 61876203 and 11801463) and the Fundamental Research Funds for the Central Universities (Nos. ZYGX2016J132 and JBK1902028).